\documentclass{amsart}
\usepackage{graphicx}
\usepackage{psfrag}
\usepackage{amsfonts}
\usepackage{amssymb}
\usepackage{amsmath}

\newtheorem{theorem}{Theorem}
\theoremstyle{plain}
\newtheorem{corollary}{Corollary}
\newtheorem{definition}{Definition}

\newtheorem{lemma}{Lemma}
\newtheorem{proposition}{Proposition}
\numberwithin{equation}{section}

\DeclareMathOperator{\length}{length}

\begin{document}
\title{Minimal Paths on Some Simple Surfaces with Singularities}
\author{Joel B. Mohler}
\address{Department of Mathematics, Lehigh University, Bethlehem, PA 18015}
\email{jbm5@lehigh.edu}
\urladdr{http://www.lehigh.edu/\symbol{126}jbm5/}
\author{Ron Umble}
\address{Department of Mathematics, Millersville University, Millersville, PA 17551}
\email{ron.umble@millersville.edu}
\urladdr{http://www.millersville.edu/\symbol{126}rumble}
\thanks{This paper is in final form and no version of it will be submitted for
publication elsewhere.}
\date{January 15, 2007}
\subjclass[2000]{ 53A04, 53A05}
\keywords{minimal path, geodesic.}

\begin{abstract}
Given two points on a soup can or conical cup with lid, we find and classify
all paths of minimal length connecting them. When the number of minimal
paths is finite, there are at most four on a can and three on a cup. At
worst, minimal paths are piece-wise smooth with three components, each of
which is a classical geodesic.  Minimal paths are geodesics in the sense of
Banchoff \cite{Banchoff}.
\end{abstract}

\maketitle

\section{Introduction}

This paper considers the following problem: {\em Given two points $A$ 
and $B$ on a soup can or conical cup (with lid) $S$, find all paths of 
minimal length connecting them.} The fact that such a path exists follows 
from Ascoli's Theorem (e.g. \cite{Munkres}):  Let $h$ be the height 
and $d$ the diameter of $S$. Let $M=2\left(h+d\right)$ and consider the 
set $G$ of all piecewise smooth constant speed paths 
$\gamma :[0,1]\rightarrow S$ such that $\gamma \left( 0\right) =A$, 
$\gamma \left( 1\right) =B$, $\length\left(\gamma \right) \leq M$ and 
$\left\Vert \gamma ^{\prime }\right\Vert \leq M$; then the speed is 
uniformly bounded and $G$ is equicontinuous.  Since equicontinuity 
is preserved in the closure, $\overline{G}$ is compact and contains a 
minimizer.

In this paper we find and classify all minimal paths from $A$ to $B$. When
the number of minimal paths is finite, there are at most four on a can and
three on a cup. At worst, minimal paths are piecewise smooth with three
components, each of which is a classical geodesic. We adapt Banchoff's
notion of a geodesic \cite{Banchoff} and prove a Hopf-Rinow-like theorem
assuring us that minimal paths are always geodesics (e.g. \cite{ONeill}).
While soup cans and conical cups are piecewise-smooth closed compact
surfaces with singularities along their rims, they have some interesting
distinguishing properties. For example, $K$ is the cone point of a conical
cup if and only if $K$ is an endpoint of every minimal path containing it.
No point on a soup can has this property.

This research grew out of an undergraduate seminar directed by the second
author in the spring of 2001 during which participants Heather Armstrong
(Cornell University), Robert Painter (College of William and Mary), Ellen
Panofsky (Lehigh University) and the first author solved the problem for
soup cans. We acknowledge the substantive contribution of each participant
to this work with thanks. Subsequently, new independent proofs of the
seminar results were obtained by the first author and appeared in his
undergraduate honors thesis \cite{Mohler} supervised by the second author.
This paper generalizes these results to conical cups in the sense that most
statements about geodesics on a conical cup have analogues on the soup can.
These analogies are realized in the limit as cone height goes to infinity.
This research was enabled by the kind and generous support of Professor
Frank Morgan, Williams College, who posed the soup can problem to our
seminar group and offered many helpful suggestions that improved this paper;
we thank him as well.

\begin{figure}[ptb]
\begin{center}
\psfrag{f}{$\phi$} \psfrag{s}{$s$} \psfrag{2pis}{$\frac{2\pi}{s}$} %
\includegraphics{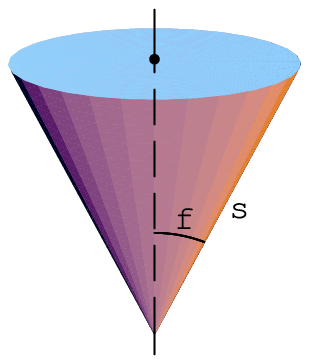}\includegraphics{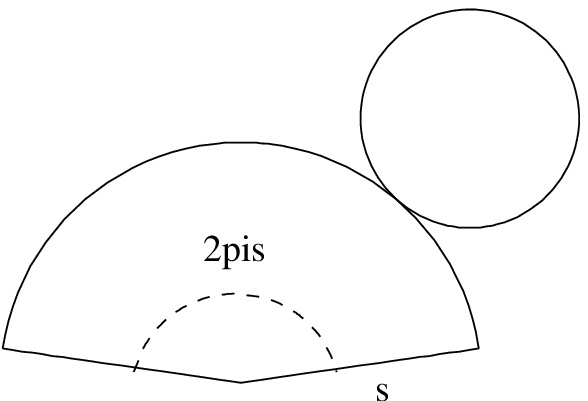}
\end{center}
\caption{A flat model of a conical cup.}
\label{fig_cone_flat}
\end{figure}

We think of a soup can or conical cup as either a surface of revolution or
an identification space formed by isometrically bending and gluing the
appropriate plane regions together along their boundaries. The
identification space point-of-view is particularly advantageous because it
provides an isometric decomposition of our surface $S$ into a family $%
\mathfrak{F}$ of compact simply-connected plane regions with piecewise
smooth boundaries; such decompositions come equipped with an arc
length-preserving quotient map $q:\mathfrak{F}\rightarrow S$. We refer to
such a pair $\left( \mathfrak{F},q\right) $ as a {\em flat model }of $S$.
In this paper, the family $\mathfrak{F}$ in a flat model of a conical cup
with cone angle $\phi $ and slant height $s=\csc \phi $ consists of a closed
unit disk $U$ (the {\em lid}) tangent to a circular sector $C$ of radius $%
s$ and subtended angle $2\pi /s$ (the {\em side}); the circle $T=\partial
U$ is called the {\em rim }(see Figure \ref{fig_cone_flat}). The family $%
\mathfrak{F}$ in a flat model of a soup can with height $h$ consists of a $%
2\pi \times h$ rectangular region $R$ (the {\em side}) and two closed
unit disks $U_{1}$ (the {\em lid}) and $U_{2}$ (the {\em base})
tangent to $R$ along its opposite edges of length $2\pi $; the circles $%
T_{i}=\partial U_{i}$ are called the {\em rims}. When there is no
confusion, we use the same symbol $X$ to denote a subset of $S$ and its
corresponding subset in the family $\mathfrak{F}$.

Let $S$ be a surface and let $I\subseteq \mathbb{R}$ be an interval. A 
{\em path} on $S$ is an element $\gamma \in C(I;S)=\{\gamma :I\rightarrow
S\mbox{ $|$ }\gamma \mbox{ continuous}\}$; in particular, if $I=\left[ a,b%
\right] $, $A=\gamma \left( a\right) $ and $B=\gamma \left( b\right) $, then 
$\gamma $ is a {\em path from $A$ to $B$}. A path $\gamma $
from $A$ to $B$ is {\em minimal} if $\length\left( \gamma
\right) \leq \length\left( \beta \right) $ for all $\beta \in C(I;S)
$. Geodesics on $S$ appear as a union of straight line segments in some
family $\mathfrak{F}$ of a flat model. In \cite{Banchoff}, Banchoff uses
flat models to study closed geodesics on \textquotedblleft
bicylinders,\textquotedblright\ but he does not consider their distance
minimizing properties. We adapt his notion of a geodesic:

\begin{definition}
\label{one}A path $\gamma \in C(\left[ a,b\right] ;S)$ is a {\em geodesic}
 if for each $t\in \left( a,b\right) $ there is a flat
model $\left( \mathfrak{F}_{t},q_{t}\right) $ of $S$ such that for every $%
\epsilon >0$ there is a $\delta >0$ and a line segment $\ell \subset 
\mathfrak{F}_{t}$ of length $\delta $ and midpoint $q_{t}^{-1}\left( \gamma
\left( t\right) \right) $ such that $\ell \subset q_{t}^{-1}\left[ \gamma
\left( t+\epsilon ,t-\epsilon \right) \right] $.
\end{definition}

\noindent We mention that Cotton et al. \cite{Cotton} use flat models to
study the classical isoperimetric problem on soup cans. The graphics in this
paper were generated by Mathematica (www.wolfram.com) and refined using XFig
(www.xfig.org).

\section{Computing minimal paths}

In this section we introduce constructive methods for finding minimal paths
and computing their length; this reduces the problem to numerical
computation. Unless explicitly stated otherwise, $S$ denotes a soup can or a
conical cup; $A$ and $B$ denote distinct points on $S$. We say that $A$ and $%
B$ are {\em separated} by $T=\partial U$ if $A\in S-U$ and $B\in U-T$ (or
vice versa). We denote the Euclidean distance between points $X$ and $Y$ by $%
XY$; if $E$ is a set, $\#E$ denotes its cardinality.

We begin with a proposition that allows us to restrict our search for
minimal paths to paths that intersect a rim at most twice:

\begin{proposition}
\label{thm_rim_crossing}Let $T$ be a rim of $S$. A path $\gamma$ on $S$ from 
$A$ to $B$ can be shortened if either

\begin{enumerate}
\item[a.] $\#\left(\gamma\cap T\right)\geq3$ or

\item[b.] $A$ and $B$ are separated by $T$ and $\#\left(\gamma\cap T\right)
\geq2$.
\end{enumerate}
\end{proposition}

\begin{proof}
(a) Consider a path $\gamma$ from $A$ to $B$.  
Suppose there are distinct points 
$X,Y,Z\in\gamma\cap T$ ordered by increasing parameter. Let $\gamma_{1}$ and
$\gamma_{2}$ be the respective pieces of $\gamma$ from $X$ to $Y$ and from $Y$
to $Z$. Then  $\triangle XYZ$ is non-degenerate and $\length\left(\gamma
_{1}\cup\gamma_{2}\right)\geq XY+YZ>XZ$, by the triangle inequality. Thus
$\left[\gamma-\left(\gamma_{1}\cup\gamma_{2}\right)\right]
\cup\overline{XZ}\ $is shorter than $\gamma$.

(b) Suppose $A$ and $B$ are separated by $T$.  With-out loss of generality, 
assume that $B\in U-T$.  By assumption, there exist distinct points 
$X,Y\in\gamma\cap T$.  Let $\gamma_{1}$ be
the piece of $\gamma$ from $X$ to $Y$ and let $\gamma_{2}$ be the piece of 
$\gamma$ from $Y$ to $B$.  Then $\triangle XYB$ is non-degenerate or $B$ 
lies in the interior of $\overline{XY}$.  In either case, 
$\length\left(\gamma_{1}\cup\gamma_{2}\right)= XY+YB>XB$ by the triangle
inequality and $\left(\gamma-(\gamma_{1}\cup\gamma_{2})\right)\cup\overline
{XB}$ is shorter than $\gamma$.
\end{proof}

\begin{corollary}
\label{endpoints}If $A$ and $B$ are separated by $T$, every minimal path
from $A$ to $B$ intersects $T$ exactly once.
\end{corollary}

\begin{proof}
If $A$ and $B$ are separated by $T$, then $\#\left(\gamma\cap T\right)
\geq1$ by continuity of $\gamma$.
\end{proof}

It is essential to know where to look for minimal paths. {\em Axial points%
} of $S$ lie on its axis of revolution. If $P$ is a non-axial point of $S$,
let $\Gamma_{P}$ denote the closed half-plane containing $P$ and bounded by
the axis.

\begin{proposition}
\label{lem_same_angle}If $A,B\in\Gamma_{Q}$, then $\Gamma_{Q}\cap S$
contains a minimal path $\gamma$ from $A$ to $B$. Furthermore, if $A$ or $B$
is non-axial, then $\gamma$ is the unique minimal path from $A$ to $B$.
\end{proposition}

\begin{proof}
We prove the result for conical cups; the proof for soup cans is similar and
left to the reader. If $A$ and $B$ lie on the lid $U$, then $\overline{AB}$ is
contained in the radius $\Gamma_{Q}\cap U$ and is (uniquely) minimal. If $A$
and $B$ lie on the side $C$, then $\overline{AB}$ is contained in the ruling
$\Gamma_{Q}\cap C$ and is (uniquely) minimal. If $A$ and $B$ are separated by
$T$ with $B\in U-T$, let $X=\Gamma_{Q}\cap T$ and let $\gamma=\overline
{AX}\cup\overline{XB};$ we claim that $\gamma$ is minimal. Let $\alpha$ be any
path from $A$ to $B;$ then $\alpha$ intersects $T$ at least once by
Proposition \ref{thm_rim_crossing}. Let $\left(\mathfrak{F},q\right)$ be
the flat model in which the side and lid are tangent at $X$ and let
$Y\in\alpha\cap T$. Then $\overline{AY}$ and $\overline{YB}$ correspond to
classical distance minimizing geodesics on $C$ and $U$, respectively. Let
$\alpha_{1}$ and $\alpha_{2}$ be the respective pieces of $\alpha$ from $A$ to
$Y$ and from $Y$ to $B$, then $\length\left(\alpha_{1}\right)\geq AY$ and
$\length\left(\alpha_{2}\right)\geq YB$. But $AY\geq AX$ since
$\overline{AX}$ is normal to the arc of sector $C;$ and $YB\geq XB$ since
$\overline{XB}$ is normal to $T$. Therefore $\length\left(\alpha\right)\geq
AY+YB\geq AB=\length\left(\gamma\right)$ and $\gamma$ is minimal.
Furthermore, if $A$ is non-axial and $Y\neq X$, segment $\overline{AY}$ is not
normal to the arc of sector $C$, in which case $AY>AX$ and $\gamma$ is
unique.
\end{proof}

Let $A$ and $B$ be non-axial points of $S$ and let $\theta$ be the angle
subtended by $\Gamma_{A}$ and $\Gamma_{B}$. If $0<\theta<\pi$, let $%
\Omega_{A,B}$ denote the closed region of space bounded by $%
\Gamma_{A}\cup\Gamma_{B}$ and subtending an angle $\theta$. If $\theta=0$,
define $\Omega_{A,B}=\Gamma_{A}=\Gamma_{B}$. If $\theta=\pi$, define $%
\Omega_{A,B}$ to be either (arbitrarily chosen) closed half-space bounded by
the plane $\Gamma_{A}\cup$ $\Gamma_{B}$. If exactly one of $A$ or $B$ is
non-axial, call it $Q$ and define $\Omega_{A,B}=\Gamma_{Q}$. If both $A$ and 
$B$ are axial, let $Q$ be an arbitrarily chosen non-axial point of $S$ and
define $\Omega_{A,B}=\Gamma_{Q}$. Below we observe that a minimal path from $%
A$ to $B$ lies in some set $\Omega_{A,B}$. This fact is well-known for
classical distance minimizing geodesics on the lid or side of a lid-less
cone and is established for $A,B\in\Gamma_{Q}$ by Proposition \ref%
{lem_same_angle}. To this end, we first establish a Hopf-Rinow-like
connection between minimal paths and geodesics (see Theorem \ref{geodesic}
below) by appealing to some well-known properties of roulettes.

\begin{definition}
\label{roulette} Let $B\in \mathbb{R}^{2}$ be fixed with respect to a closed
convex curve $C_{1}$; let $C_{2}$ be any plane curve.  The {\em roulette 
generated by $B$} is the curve traced out by $B$ as $C_{1}$
rolls without slipping along $C_{2}$.%
\end{definition}

\begin{figure}[tbh]
\begin{center}
\psfrag{A}{$A$} \psfrag{B}{$B$} \psfrag{P}{$P$} %
\includegraphics{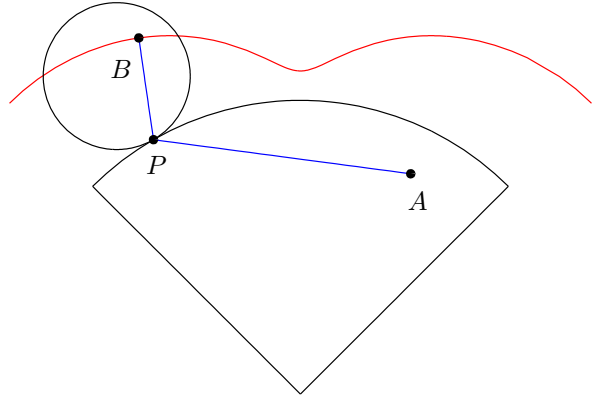}
\end{center}
\caption{A non-minimal path crossing the rim.}
\label{figure_2}
\end{figure}

\noindent In our considerations below, the point $B$ in Definition \ref%
{roulette} lies on a lid $U$ in some flat model $\left(\mathfrak{F}
_{0},q_{0}\right)$, the curve $C_{1}$ is the rim $T=\partial U$ and $C_{2}$
is an edge of the side. When $S$ is a soup can, the roulette generated by $B$
is a {\em cycloid}; when $S$ is a conical cup, the roulette generated by $%
B$ is an {\em epicycloid} (see Figure \ref{figure_2}). 

\begin{lemma}
\label{non-geodesic}If $A$ and $B$ are separated by a rim, every minimal
path from $A$ to $B$ is a geodesic.
\end{lemma}

\begin{figure}[tbh]
\begin{center}
\psfrag{A}{$A$} \psfrag{B}{$B$} \psfrag{E}{$E$} \psfrag{P}{$P$} %
\psfrag{Pp}{$P^\prime$} \psfrag{Ppp}{$P^{\prime\prime}$} \psfrag{Q}{$Q$} %
\includegraphics{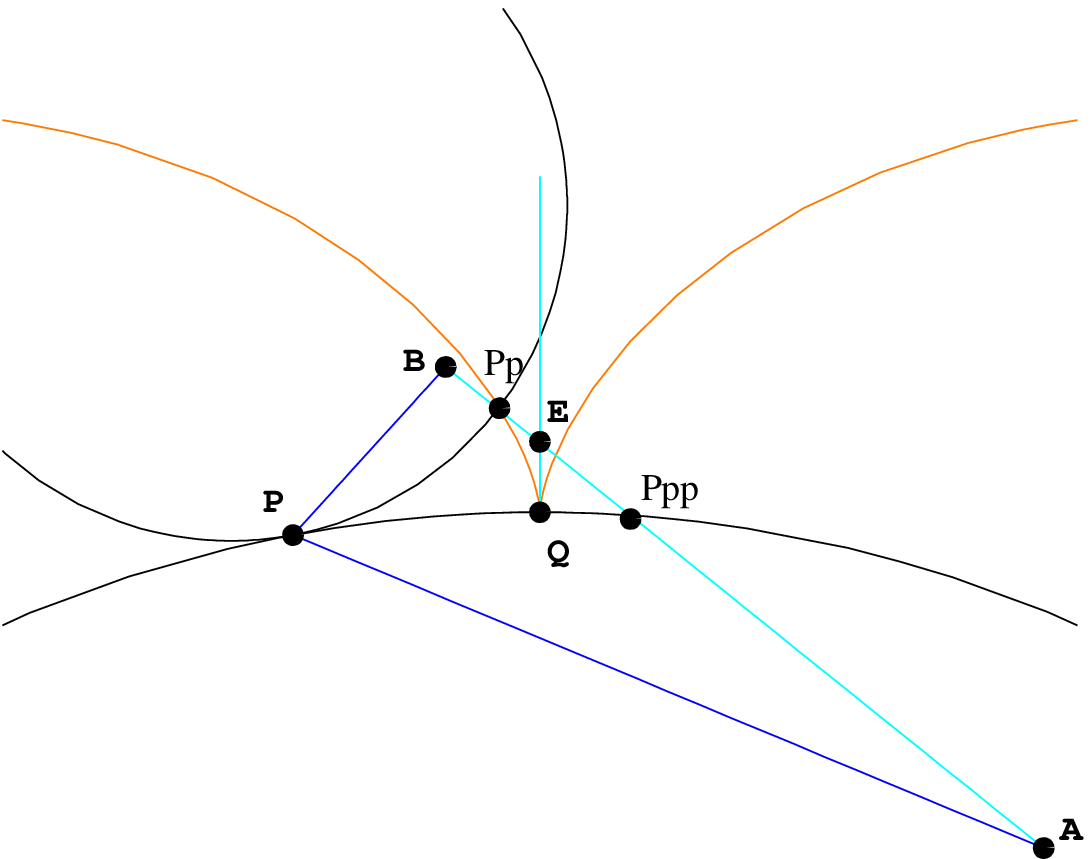}
\end{center}
\caption{Shortening a non-geodesic path from $A$ and $B$.}
\label{figure_3}
\end{figure}

\begin{proof}
We prove the result for conical cups; the result for soup cans is obtained in
the limit as cone height goes to infinity. Given a non-geodesic path $\gamma$
from $A$ to $B$ crossing the rim $T$ at point $P$, we construct an alternate
path $\gamma^{\prime}$ such that $\length\left(\gamma^{\prime}\right)
<\length\left(\gamma\right)$. Let $\gamma_{1}$ the piece of $\gamma$ from
$A$ to $P$ and let $\gamma_{2}=\gamma-\gamma_{1}$. Since $P$ is shared by the
side and lid, there are classical distance minimizing geodesics from $A$ to
$P$ and from $P$ to $B$. So assume that paths $\gamma_{i}$ are such geodesics,
in which case $\gamma$ intersects $T$ exactly once at $P$. Consider the flat
model $\left(\mathfrak{F},q\right)$ in which the lid $U$ is tangent to the
sector $C$ at $P$. Since $\gamma$ is not a geodesic, $q^{-1}\left(\gamma\right)
=\overline{AP}\cup\overline{PB}$ is not straight (see Figure \ref{figure_2}).
Construct a shorter path $\gamma^{\prime}$ from $A$ to $B$ as follows: Let
$R$ denote the arc in the boundary of sector $C$ and $T=\partial U$. Let $P^{\prime}=\overline{AB}\cap T$ and
$P^{\prime\prime}=\overline{AB}\cap R$; let $s_{1}$ be the length of the arc
along $T$ from $P$ to $P^{\prime}$ and inside $\triangle APB$; and let $s_{2}$
be the length of the arc along $R$ from $P$ to $P^{\prime\prime}$ (see Figure
\ref{figure_3}). If $s_{1}\leq s_{2}$, let $\alpha$ be the epicycloid
generated by $P^{\prime}$ as disk $U$ rolls along $R$ towards $P^{\prime\prime}$.
Let $Q$ be the first cusp that appears. We claim that $P^{\prime}%
P^{\prime\prime}>QP^{\prime\prime}$. Let $E$ be the intersection of
$\overline{P^{\prime}P^{\prime\prime}}$ and the line normal to $R$ at
$Q$. Construct the chord $\overline{QP^{\prime\prime}}$; note that $P^{\prime
}P^{\prime\prime}>EP^{\prime\prime}$ since $P^{\prime}$ and $P^{\prime\prime}$
are on opposite sides of $\overline{QE}$. Furthermore, $\theta=m\angle
EQP^{\prime\prime}>\pi/2$ so that $\cos\theta<0$. Then by the Law of Cosines,
\[
(EP^{\prime\prime})^{2}=(QP^{\prime\prime})^{2}+\left(QE\right)
^{2}+\left\vert 2(QP^{\prime\prime})\left(QE\right)\cos\theta\right\vert >(QP^{\prime\prime})^{2}%
\]
so that $EP^{\prime\prime}>QP^{\prime\prime}$, which verifies the claim. Thus
$AB>AP^{\prime\prime}+QP^{\prime\prime}+P^{\prime}B$ and $\gamma^{\prime
}=\overline{AP^{\prime\prime}}\cup\overline{P^{\prime\prime}Q}\cup
\overline{P^{\prime}B}$ is shorter than $\gamma$. On the other hand, if
$s_{1}>s_{2}$ consider the epicycloid generated by $P^{\prime\prime}$ as sector $C$
rolls along $T$ towards $P^{\prime}$ and construct the path $\gamma^{\prime
}=\overline{AP^{\prime\prime}}\cup\overline{P^{\prime}Q}\cup\overline
{P^{\prime}B}$.
\end{proof}

\begin{theorem}
\label{geodesic}Every minimal path from $A$ to $B$ is a geodesic.
\end{theorem}

\begin{proof}
Let $\gamma$ be a minimal path from $A$ to $B$. If $\gamma\subset S-U$ or
$\gamma\subset U$, the result is classical. The case with $A$ and $B$
separated by $T$ is given by Lemma \ref{non-geodesic}. Suppose that $\gamma$
intersects $T$ exactly twice at $P_{1}$ and $P_{2}$. If $A,B\in T$, then
$\gamma$ is a geodesic in $U$. So assume that $A$ lies off $T$, then $A\in
S-U$, for otherwise $\triangle AP_{1}P_{2}\subset U$ is non-degenerate and
$\gamma$ is not minimal. Choose a point $Q\in\gamma\cap\left(U-T\right)$
and let $\gamma_{1}$, $\gamma_{2}$ and $\gamma_{3}$ be the respective pieces
of $\gamma$ from $A$ to $Q$, from $Q$ to $B$ and from $P_{1}$ to $P_{2}$. Then
$\gamma_{1}$ is a geodesic by Lemma \ref{non-geodesic} and $\gamma_{3}$ is a
geodesic in $U$. Therefore $\gamma_{1}\cap\gamma_{3}=\overline{P_{1}Q}$ is a
geodesic. If $B\in T$, then $\gamma_{2}\subset\gamma_{3}$ and we're done. If
$B$ is off $T$, then $B\in S-U$ and $\gamma_{2}$ is a geodesic by Lemma
\ref{non-geodesic}. Therefore $\gamma_{2}\cap\gamma_{3}=\overline{QP_{2}}$ is
a geodesic and the conclusion follows.
\end{proof}

\begin{proposition}
\label{twelve}Let $P$ be a non-axial point of $S$. If a minimal path $\gamma$
from $A$ to $B$ intersects $\Gamma_{P}$, either $\gamma\cap\Gamma_{P}=\gamma$
or $\gamma\cap\Gamma_{P}$ is a single point.
\end{proposition}

\begin{proof}
Suppose distinct points $Q,Q^{\prime}\in\gamma\cap\Gamma_{P};$ let $\left(
\mathfrak{F},q\right)$ be the flat model containing $\Gamma_{P}\cap S$.
Since $\gamma$ is minimal, $\overline{QQ^{\prime}}\subset\gamma\cap\Gamma_{P}$
by Proposition \ref{lem_same_angle}. So consider the connected component
$\alpha$ of $\gamma\cap\Gamma_{P}\subset\mathfrak{F}$ containing
$\overline{QQ^{\prime}};$ let $A=\gamma\left(a\right)$, $B=\gamma\left(
b\right)$. Since $\gamma$ is a geodesic by Theorem \ref{geodesic}, no point
$E=\gamma\left(t\right)$, $t\in\left(a,b\right)$ is an endpoint of
$\alpha$, for if it were, every neighborhood of $E$ would violate the
requirements of Definition \ref{one}. Therefore $A$ and $B$ are the endpoints
of $\alpha$ and $\gamma=\alpha$.
\end{proof}

\begin{theorem}
\label{thm_interior_angle}Every minimal path from $A$ to $B\ $lies in some
set $\Omega_{A,B}$.
\end{theorem}

\begin{figure}[tbh]
\begin{center}
\psfrag{Pi}{$\Pi$} \includegraphics{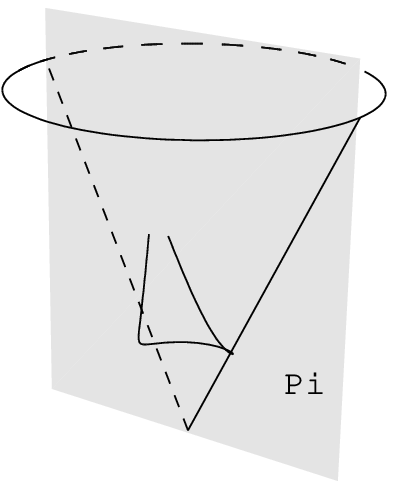}\hspace{.25in}%
\includegraphics{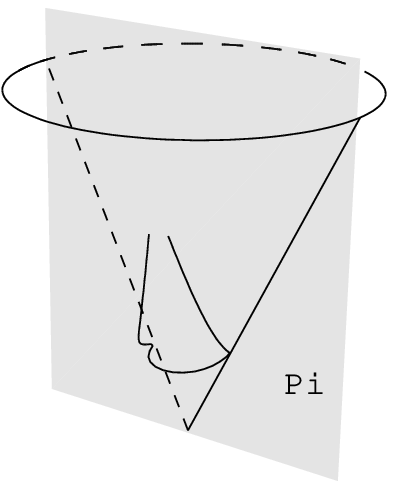}
\end{center}
\caption{A geodesic reflected in vertical plane $\Pi$.}
\label{fig_geo_reflected}
\end{figure}

\begin{proof}
The case of $A,B\in\Gamma_{Q}$ for some $Q$ was established in Proposition
\ref{lem_same_angle}. So assume $A$ and $B$ are non-axial points of $S$ and
let $\theta$ be the angle between $\Gamma_{A}$ and $\Gamma_{B}$. If
$0<\theta<\pi$, suppose $\gamma\cap\Omega_{A,B}^{c}\neq\varnothing$. Then by
Proposition \ref{twelve}, we may assume that $\gamma\subset\left\{
A,B\right\}  \cup\Omega_{A,B}^{c}$. Since $0<\theta<\pi$, there is a plane
$\Pi$ containing the axis and bounding an open half-space $H$\ with $A,B\in
H$. Let $\overline{\gamma}$ be the reflection of $\gamma\cap H^{c}$ in $\Pi;$
then $\beta=\overline{\gamma}\cup\left(\gamma\cap H\right)$ is a new path
passing through $\Gamma_{A}$ (and $\Gamma_{B}$) at least twice, i.e., at
$A=\beta\left(a\right)$ and $P=\beta\left(u\right)$ for some $u\neq a$
(see Figure \ref{fig_geo_reflected}). If $A=P$, we can shorten $\beta$ by
redefining $\beta\left(t\right)=A$ for all $t\in\left[a,u\right];$
otherwise we can shorten $\beta$ by applying Proposition \ref{lem_same_angle},
contradicting minimality in either case. Finally, suppose that $\theta=\pi$
and $\gamma$ contains points on both sides of plane $\Pi=\Gamma_{A}\cup
\Gamma_{B}$. If $A=\gamma\left(a\right)$ and $\gamma\left(t\right)
\in\Gamma_{A}$ for some $t\neq a$, then $\gamma\subset\Gamma_{A}$ by
Proposition \ref{twelve}, which is a contradiction. Therefore $\gamma
\cap\Gamma_{A}=A$. Likewise, $\gamma\cap\Gamma_{B}=B$ and the conclusion
follows by the continuity of $\gamma$.
\end{proof}

\noindent Note that when $\theta=\pi$ and a half-space $\Omega_{A,B}$ has
been chosen, $\Omega_{A,B}$ and the closure of its complement each contain a
minimal path from $A$ to $B$. With the characterization of minimal length
paths developed to this point, we now exhibit numerical methods to find such
paths in Propositions \ref{thm_cone_lid_side}-\ref{thm_lid_base}. Denote the
cylindrical coordinates of point $P$ by $P\left[r,\theta_{P},z\right]$.

\begin{proposition}
\label{thm_cone_lid_side}Let $A$ and $B$ be respective points on the side
and lid of $S$. If $A$ and $B$ are axial points, there is a minimal path
from $A$ to $B$ in $\Gamma_{Q}$ for each non-axial $Q$. If exactly one of $A$
or $B$ is non-axial, call it $Q,$ there is a unique minimal path from $A$ to 
$B$ in $\Gamma_{Q}$. If $A$ and $B$ are non-axial and separated by rim $T$,
position $S$ so that $A$ has cylindrical angle $\theta_{A}=0$, let $h$ be
the height of $S$, let $b$ be the distance from $B$ to the axis, let $%
\theta\in\left[ 0,\pi\right]$ be the angle between $\Gamma_{A}$ and $%
\Gamma_{B}$ and let $t\in\left[0,\theta\right];$ let $\gamma_{1}$ be the
classical distance minimizing geodesic from $A$ to $P\left[1,t,h\right]$ and
let $\gamma _{2}=\overline{PB}$.

\begin{enumerate}
\item[a.] If $S$ is a conical cup with slant height $s$ and $a$ is the
distance from $A$ to the cone point, then%
\begin{align*}
\ell(t) & =\length\left(\gamma_{1}\cup\gamma_{2}\right) \\
& =\sqrt{a^{2}-2as\cos\left(t/s\right)+s^{2}}+\sqrt{b^{2}-2b\cos (\theta-t)+1%
}.
\end{align*}

\item[b.] If $S$ is a soup can and $a$ is the distance from $A$ to the rim,
then%
\begin{align*}
\ell\left(t\right) & =\length\left(\gamma_{1}\cup\gamma_{2}\right) \\
& =\sqrt{a^{2}+t^{2}}+\sqrt{b^{2}-2b\cos(\theta-t)+1}.
\end{align*}
\end{enumerate}

\noindent In either case, $\ell$ is minimized at some $t_{1}\in\left[
0,\theta\right]$ and a path $\gamma_{1}\cup\gamma_{2}$ of length 
$\ell\left(t_{1}\right)$ is a minimizer.
\end{proposition}

\begin{proof}
The axial cases were established by Proposition \ref{lem_same_angle}; so
assume $A$ and $B$ are non-axial. The formulation for length $\ell\left(
t\right)$ follows from the Law of Cosines; the restriction of parameter
$t\in\left[0,\theta\right]$ follows from Theorem \ref{thm_interior_angle}.
The conclusion follows from continuity of $\ell$ on $\left[0,\theta\right]
$.
\end{proof}

\noindent Our next result is specific to conical cups:

\begin{proposition}
\label{cone_side}Suppose $A$ and $B$ are non-axial points on the side of a
conical cup $S$ with height $h$ and slant height $s$, positioned so that $A$
has cylindrical angle $\theta_{A}=0$. Let $\theta\in\left[0,\pi\right]$ be
the angle between $\Gamma_{A}$ and $\Gamma_{B}$ and let $t,u\in\left[
0,\theta\right];$ let $\gamma_{1}$ be the classical minimal geodesic from $A$
to the point $P_{1}\left[1,t,h\right]\in T$, let $\gamma_{2}$ be the
classical minimal geodesic from $P_{2}\left[1,u,h\right]\in T$ to $B$ and
let $\gamma_{3}=\overline{P_{1}P_{2}}$ on the lid. If the respective
distances from $A$ and $B$ to the cone point are $a$ and $b$, then 
\begin{align*}
\ell(t,u) & =\length\left(\gamma_{1}\cup\gamma_{2}\cup\gamma_{3}\right) \\
& =\sqrt{a^{2}-2as\cos(t/s)+s^{2}}+\sqrt{b^{2}-2bs\cos\left(u/s\right) +s^{2}%
} \\
& \hspace*{0.75in}+\sqrt{2-2\cos\left(\theta-t-u\right)}.
\end{align*}
Candidates for minimizer appear as:

\begin{enumerate}
\item[a.] a classical minimal geodesic $\alpha$ on the side, or

\item[b.] a path $\beta=\gamma_{1}\cup\gamma_{2}\cup\gamma_{3}$ across the
lid of minimal length $\ell\left(t_{1},u_{1}\right)$.

\item[ ] A minimal path from $A$ to $B$ is minimal in $\left\{ \alpha
,\beta\right\} $.
\end{enumerate}
\end{proposition}

\begin{proof}
The result follows by continuity of $\ell$ on the compact set $\left[
0,\theta\right]\times\left[0,\theta\right]$.
\end{proof}

\noindent Note that one or both of $A$ and $B$ in Proposition \ref{cone_side}
may lie on the rim. This result generalizes to soup cans by applying the
same analysis at both ends:

\begin{proposition}
\label{can_side}Suppose $A$ and $B$ are points on the side of a soup can $S$
with height $h$, positioned so that $A$ has cylindrical angle $\theta_{A}=0$%
. Let $\theta\in\left[0,\pi\right]$ be the angle between $\Gamma_{A}$ and $%
\Gamma_{B}$ and let $t,u\in\left[0,\theta\right];$ let $\gamma_{1}$ be the
classical minimal geodesic from $A$ to the point $P_{1}\left[ 1,t,h\right]%
\in T$, let $\gamma_{2}$ be the classical minimal geodesic from $P_{2}\left[%
1,u,h\right]\in T_{1}$ to $B$ and let $\gamma_{3}=\overline{P_{1}P_{2}}$ on
the lid. If $a$ and $b$ are the respective distances from $A$ and $B$ to rim 
$T_{1}$. Then%
\begin{align*}
\ell_{1}\left(t,u\right) &
=\length\left(\gamma_{1}\cup\gamma_{2}\cup\gamma_{3}\right) \\
& =\sqrt{a^{2}+t^{2}}+\sqrt{b^{2}+u^{2}}+\sqrt{2-2\cos(\theta-t-u)}.
\end{align*}
Furthermore, let $\gamma_{4}$ be the classical minimal geodesic from $A$ to $%
Q_{1}\left[1,t,0\right]\in T_{2}$, let $\gamma_{5}$ be the classical minimal
geodesic from $Q_{2}\left[1,u,0\right]\in T_{2}$ to $B$ and let $\gamma_{6}=%
\overline{Q_{1}Q_{2}}$ on the base. Then%
\begin{align*}
\ell_{2}\left(t,u\right) &
=\length\left(\gamma_{4}\cup\gamma_{5}\cup\gamma_{6}\right) \\
& =\sqrt{\left(h-a\right)^{2}+t^{2}}+\sqrt{\left(h-b\right)^{2}+u^{2}} \\
& \hspace*{0.75in}+\sqrt{2-2\cos(\theta-t-u)}.
\end{align*}
Candidates for minimizer appear as:

\begin{enumerate}
\item[a.] a classical minimal geodesic $\alpha$ on the side,

\item[b.] a path $\beta_{1}=\gamma_{1}\cup\gamma_{2}\cup\gamma_{3}$ across
the lid of minimal length $\ell_{1}\left(t_{1},u_{1}\right)$, or

\item[c.] a path $\beta_{2}=\gamma_{4}\cup\gamma_{5}\cup\gamma_{6}$ across
the base of minimal length $\ell_{2}\left(t_{2},u_{2}\right)$.

\item[ ] A minimal path from $A$ to $B$ is minimal in $\left\{ \alpha
,\beta_{1},\beta_{2}\right\} $.
\end{enumerate}
\end{proposition}

\begin{proof}
The result follows by continuity of $\ell_{i}$ on the compact set $\left[
0,\theta\right]\times\left[0,\theta\right]$.
\end{proof}

\noindent Again, one or both of $A$ and $B$ in Proposition \ref{can_side}
may lie on a rim. Our next result is specific to soup cans.

\begin{proposition}
\label{thm_lid_base} Let $A$ and $B$ be non-axial points on the lid and base
of a soup can $S$ whose respective distances from the axis are $a$ and $b$.
Position $S$ so that $A$ has cylindrical angle $\theta_{A}=0$. Let $h$ be
the of height of $S$, let $\theta$ be the angle between $\Gamma_{A}$ and $%
\Gamma_{B}$, let $t,u\in\left[0,\theta\right]$, and consider the point $P%
\left[1,t,h\right]\in T_{1}$ and $Q\left[1,u,0\right]\in T_{2}$. Let $%
\gamma_{1}=\overline{AP}$ on the lid, let $\gamma_{2}=\overline{QB}$ on the
base and let $\gamma_{3}$ be a classical minimal geodesic on the side from $P
$ to $Q$. Then 
\begin{align*}
\ell\left(t,u\right) & =\length\left(\gamma_{1}\cup\gamma_{2}\cup
\gamma_{3}\right) \\
& =\sqrt{a^{2}-2a\cos(\theta-t)+1}+\sqrt{b^{2}-2b\cos(\theta-u)+1} \\
& \hspace*{0.75in}+\sqrt{h^{2}+(\theta-t-u)^{2}}.
\end{align*}
Then $\ell$ is minimized at some $\left(t_{1},u_{1}\right)\in\left[ 0,\theta%
\right]\times\left[0,\theta\right]$ and a path $\gamma_{1}\cup\gamma_{2}\cup%
\gamma_{3}$ of length $\ell\left(t_{1},u_{1}\right)$ is a minimizer.
\end{proposition}

\begin{proof}
Again, the result follows from continuity of $\ell$.
\end{proof}

We conclude this section with an interesting observation that distinguishes
conical cups from soup cans.

\begin{definition}
Non-axial points $A$ and $B$ on a surface of revolution are {\em diaxial} 
if the angle between $\Gamma_{A}$ and $\Gamma_{B}$ is $\pi$.
\end{definition}

\begin{proposition}
\label{two}A point $K$ on a conical cup is the cone point if and only if $K$
is an endpoint of every minimal path containing it. Furthermore, every
minimal path ending at the cone point lies in some plane containing the axis.
\end{proposition}

\begin{figure}[tbh]
\begin{center}
\psfrag{pisinphi}{$\pi\sin\phi$} \psfrag{A}{$A$} \psfrag{B}{$B$} %
\psfrag{K}{$K$} \includegraphics{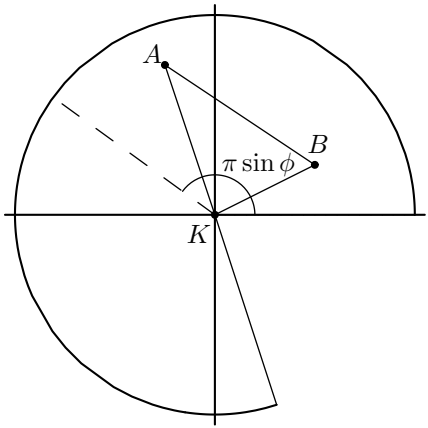}
\end{center}
\caption{A minimal path on the side misses the cone point.}
\label{fig_void_vertex}
\end{figure}

\begin{proof}
Suppose that $K$ is an endpoint of every minimal path containing it. But $K$
cannot lie on the side since some ruling would contain it; $K$ cannot lie on
the rim since the path $\Gamma_{K}\cap S$ would contain it; and $K$ cannot lie
on the lid since some diameter would contain it. Thus $K$ is the cone
point.\ Conversely, we claim that the cone point $K$ is an endpoint of every
minimal path containing it. Suppose that $\gamma$ is a minimal path from $A$
to $B$, where $A,B\neq K$. Since minimal paths between two points on the lid
or from a point on the side to a point on the lid miss $K$, we may assume that
$A$ and $B$ lie on the side. If $A$ and $B$ are diaxial, let $\alpha$ be a
ruling that misses $A$ and $B$. If not, choose a plane $\Pi$ containing the
axis and bounding an open half-space $H$ with $A,B\in H$ and let $\alpha$ be a
ruling in $\Pi\cap S$. In either case, choose a flat model $\left(
\mathfrak{F},q\right)$ with $q$ mapping the edges of sector $C$ onto
$\alpha$. If $\phi$ is the cone angle, $m\angle AKB\leq\pi\sin\phi<\pi$ so
that $\overline{AB}\subset C;$ and furthermore, $K\notin\overline{AB}$ since
$K\subset\partial C$. It follows that $\gamma=q\left(\overline{AB}\right)$
misses $K$ (see Figure \ref{fig_void_vertex}). Finally, suppose $A$ is the
cone point. If $B$ is non-axial, there is a minimal path from $A$ to $B$ in
$\Gamma_{B}\cap S$. If $B$ is axial, choose any point $P$ distinct from $A$
and $B;$ the path $\Gamma_{P}\cap S$ is minimal. Thus every minimal path
beginning at the cone point lies in some plane containing the axis.
\end{proof}

\section{Non-unique minimal paths}

\begin{figure}[tbh]
\begin{center}
\psfrag{A}{$A$} \psfrag{B}{$B$} \includegraphics{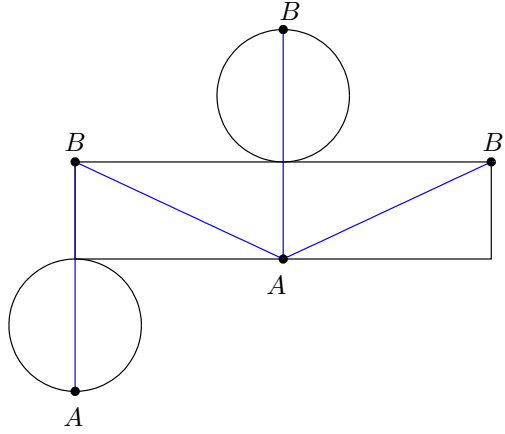}
\end{center}
\caption{Four minimal paths from $A$ to $B$.}
\label{fig_can_four_geo}
\end{figure}

We conclude the paper with a discussion of some situations in which multiple
minimal paths from $A$ to $B$ exist. Consider a soup can $S$ with diaxial
points $A$ and $B$ on opposite rims. If height $h\,<\frac{\pi^{2}-4}{4}$,
there are exactly two minimal paths from $A$ to $B$--one across the base and
one across the lid. If $h>\frac{\pi^{2}-4}{4}$, again there are exactly two
minimal paths, but this time both are on the side. And if $h=\frac{\pi^{2}-4%
}{4}$, there are exactly four--one across the lid, one across the base and
two on the side (see Figure \ref{fig_can_four_geo}). In general, given
diaxial points $A$ and $B$, there are at most four minimal paths connecting
them; in particular, if height $h>\frac{\pi^{2}-4}{4}$, there are at most
three--one across the lid and two around the side. This motivates the
following somewhat surprising result for conical cups.

\begin{theorem}
The side of a conical cup with slant height $s$ contains diaxial points $A$
and $B$ joined by three minimal paths. Specifically, let $r_{\max}$ be the
maximum radius of curvature along an epicycloid generated by a point on the
rim as it rolls along the sector arc in a flat model. If $A$ is non-axial
point of distance $a\geq s-r_{\max}+2$ from the cone point, there are three
minimal paths from $A$ to the diaxial point $B$ whose distance $b$ from the
cone point is%
\begin{equation}
b=\frac{2\left(s-a+1\right)\left(s+1\right)}{2s+2-a\left[ 1+\cos\left(\pi/s%
\right)\right]}.   \label{relation}
\end{equation}
\end{theorem}

\begin{figure}[tbh]
\begin{center}
\psfrag{A}{$A$} \psfrag{B}{$B$} \includegraphics{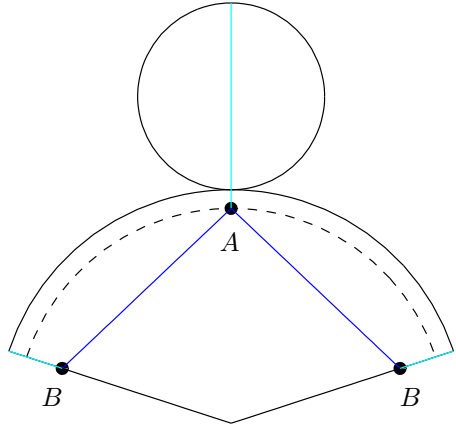}
\end{center}
\caption{Three minimal paths from $A$ to $B$.}
\label{fig_cup_three_geo}
\end{figure}

\begin{proof}
First observe that minimal paths across the lid from $A$ to $B$ follow a
diameter: Consider the flat model in which sector $C$ is centered at the
origin with one edge along the positive $x$-axis, the disk $U$ is centered at
$\left(s+1,0\right)$, $A\left(s-a,0\right)$ and $B\left(\left(
s-b\right)\cos\left(\pi/s\right),\left(s-b\right)\sin\left(
\pi/s\right)\right)$. Let $\alpha$ be the trace of the epicycloid
generated by $P(s+2,0)$ as $U$ rolls along $C$. Since the radius of curvature
along $\alpha$ attains its maximum $r_{\max}$ at $P$, the circle centered at
$A$ of radius $s-a+2\leq r_{\max}$ intersects $\alpha$ only at $P$ and the
distance from $A$ to $\alpha$ is minimized at $P$. Of course, the distance
from $B$ to the sector arc is minimized at $Q\left(s\cos\left(
\pi/s\right),s\sin\left(\pi/s\right)\right)$, so minimal paths across
the rim from $A$ to $B$ lie in $\left(\Gamma_{A}\cup\Gamma_{B}\right)\cap
S$ as claimed. To obtain relation (\ref{relation}), note that sector $C$
subtends angle $2\pi/s$ and contains $\overline{AB}$. Since $A$ and $B$ are
diaxial, $m\angle AOB=\pi/s$ and%
\[
AB=\sqrt{a^{2}-2ab\cos\left(\pi/s\right)+b^{2}}%
\]
by the Law of cosines. Now assume for the moment that $AB$ is the minimal
length $2s-a-b+2$ of paths across the lid. Then%
\[
(2s-a-b+2)^{2}=a^{2}-2ab\cos\left(\pi/s\right)+b^{2}%
\]
and solving for $b$ gives relation (\ref{relation}). So given $A$ such that
$a\geq s-r_{\max}+2$, let $b$ be given by relation (\ref{relation}). Choose
$B$ at distance $b$ from the cone point and positioned so that $A$ and $B$ are
diaxial. Then $A$ and $B$ are joined by exactly three minimal paths (see
Figure \ref{fig_cup_three_geo}).
\end{proof}

Sufficiently tall soup cans admit an analogous result.

\begin{corollary}
Let $A$ and $B$ be diaxial points on the side of a soup can $S$ whose
respective distances from the lid are $c$ and $d$. If $S$ has critical
height 
\begin{equation*}
h=\frac{\pi^{2}+4c^{2}-4}{4c+4}
\end{equation*}
and $d=h-c$, there are exactly four minimal paths from $A$ to $B$. In
particular, if $A$ and $B$ lie on opposite rims, $h=\frac{\pi^{2}-4}{4}$; if 
$A$ and $B$ lie on the mid-circle, $h=\pi-2$.
\end{corollary}

\begin{proof}
Set $a=s-c$ and $b=s-d$ in relation (\ref{relation}) and obtain
\[
d=-\frac{2c+2+s\left(c-s\right)\left[1-\cos\left(\pi/s\right)
\right]}{2s+2+\left(c-s\right)\left[1+\cos\left(\pi/s\right)
\right]}.
\]
Since%
\[
-\lim_{s\rightarrow\infty}\frac{2c+2+s\left(c-s\right)\left[
1-\cos\left(\pi/s\right)\right]}{2s+2+\left(c-s\right)\left[
1+\cos\left(\pi/s\right)\right]}=\frac{\pi^{2}-4c-4}{4c+4},
\]
the critical height $h=c+\frac{\pi^{2}-4c-4}{4c+4}$.
\end{proof}

Multiple minimal paths can join points in more general position. For
example, let $P$ and $B$ be points on the rim of a soup can $S$ such that
the angle $\theta$ between $\Gamma_{P}$ and $\Gamma_{B}$ satisfies $%
0<\theta_{P}<\pi/2$. Consider the flat model whose lid is tangent to the
rectangular side $R$ at $P$, and let $B^{\prime}\in\partial R$ be the point
identified with $B$. For $S$ with sufficient height $h$, there is a point $%
A\in\overleftrightarrow {PB}\cap R$ such that $\triangle ABB^{\prime}$ is
isosceles and there are at least two minimal paths in $\Omega_{A,B}$ from $A$
to $B$--one across the lid and one around the side. In fact, there is a
critical height $h$ at which there are three minimal paths--the two just
mentioned and a third across the base. Computing this critical height $h$
requires more machinery than we needed above, but is similar in spirit. We
conclude with a computation of $h$ when $A$ and $B$ lie on the mid-circle of
a soup can. We need the following well-known fact about roulettes:

\begin{proposition}
Given a point $B$ and curves $C_{1}$ and $C_{2}$ as in Definition \ref%
{roulette}, let $\beta$ be the roulette generated by $B$. If $\left\Vert
\beta^{\prime}\left(t\right)\right\Vert >0$, the line normal to $\beta$ at $%
\beta\left(t\right)$ passes through the point of tangency $\alpha
_{2}\left(t\right)$.
\end{proposition}

\noindent Now if $A$ and $B$ are joined by a minimal path $\gamma_{1}$
crossing the rim $T_{1}$ at $P$ and $Q$, the mirror image of $\gamma_{1}$ is
a minimal path $\gamma_{2}$ across the base $U_{2}$. We wish to compute the
critical height at which there is a third path $\gamma_{3}$ of the same
length around the side. Assume that $Q$ is the point on $T_{1}$ closest to $B
$, and let $\beta$ be the cycloid generated by $Q$ as $U_{1}$ rolls along
its edge of rectangle $R$ in a flat model. Let $n$ be the line normal to $%
\beta$ passing through $A$ and let $\theta$ be its angle of inclination.
When $T_{1}$ is tangent to $R$ at $P$, the central angle in $T_{1}$
subtended by the chord of $n$ is $2\theta$ and the height of the soup can is 
$h=2BP\sin\theta$. So if $h$ is the desired critical height, $%
\theta=\sin^{-1}\left(\frac{h}{2BP}\right)$ is a solution of%
\begin{equation}
\theta-\sin\theta=BP\left(1-\cos\theta\right).   \label{eleven}
\end{equation}
Conversely, solving equation (\ref{eleven}) for $\theta$ determines $h$.
Other cases are similar in spirit and left to the reader.

\end{document}